\documentclass{article}

\everydisplay{\mathsurround 0pt} 

\usepackage[a4paper]{geometry}
\usepackage{amsmath,amssymb,amsthm,mathtools}
\usepackage[colorlinks,citecolor=blue,urlcolor=blue]{hyperref}
\usepackage[affil-it]{authblk}
\usepackage[round]{natbib}
 \newcommand{\E}{\mathbb{E}}

\newcommand{\AAA}{\mathcal{A}}
  \newcommand{\FFF}{\mathcal{F}}

\newcommand{\R}{\mathbb{R}}

\newtheorem{lem}{Lemma}

\newtheorem{theorem}{Theorem}

\title{Gauss and the identity function - a tale of characterizations of the normal distribution} 

\author{Christophe Ley
\thanks{Electronic address: \texttt{christophe.ley@ugent.be}. I would like to thank Marc Hallin and Yvik Swan for inspiring discussions about the normal distribution.}}
\affil{Ghent University, Ghent, Belgium}

\date{}

\begin{document}

\maketitle

\begin{abstract}

\indent The normal distribution is well-known for several results that it is the only to fulfil. The aim of the present paper is to show that many of these characterizations actually follow from the fact that the derivative of the log-density of the normal distribution is the (negative) identity function. This  \emph{a priori} very simple yet surprising observation allows a deeper understanding of existing characterizations and paves the way to an immediate extension to a general density $x\mapsto p(x)$ by replacing $-x$ in these results with $(\log p(x))'$.

\end{abstract}


\section{INTRODUCTION}

The normal or Gaussian distribution is the most popular probability law in statistics and probability. The reasons for this popularity are manifold, including the nice bell curve shape, the simple form of the density
$$
x\mapsto \phi_{\mu,\sigma}(x):=\frac{1}{\sqrt{2\pi}\sigma}\exp\left(-\frac{(x-\mu)^2}{2\sigma^2}\right), \quad x\in\R,
$$
with easily interpretable location parameter $\mu\in\R$ and scale parameter $\sigma>0$, the ensuing mathematical tractability, the straightforward extension to the multivariate normal density (which we shall however not deal with in this paper) or the fact of being the limit distribution in the Central Limit Theorem. Besides these major appeals, the normal distribution is also famous for satisfying various characterizations, the latter being theoretical results that only one distribution (or one class of distributions) is fulfilling. Carl Friedrich Gauss himself has obtained the normal density by searching for a probability distribution where the maximum likelihood estimator of the location parameter \emph{always} (see Section~\ref{sec:MLE} for a precise meaning) coincides with the most intuitive estimator, namely the sample average. Numerous other characterizations of this popular distribution have followed, and in general it took the researchers decades to extend them to other distributions, often in an ad hoc way.

In the present paper, we will show that an apparently innocuous characterization of the normal distribution turns out to be a crucial building block in several more famous characterizations. This characterization is the fact that $(\log \phi_{0,1}(x))'=-x$ or, equivalently, $\frac{d}{d\mu}(\log \phi_{\mu,\sigma}(x))=\frac{x-\mu}{\sigma^2}$. In the former notation we speak of the derivative of the log-density, while  the second case features the location score function (we will refer to both settings as the ``identity function'' or ``location score function''). It is straightforward to see that the normal distribution is the only one for which these results hold.  We shall show in the remainder of this paper that this particular  characterization of the normal distribution via the identity function actually lies at the core of many characterizations that convey its special role to the normal distribution. In particular, we shall focus on the already mentioned maximum likelihood characterization (Section~\ref{sec:MLE}), on a singular Fisher information matrix characterization in skew-symmetric distributions (Section~\ref{sec:info}), on Stein characterizations (Section~\ref{sec:Stein}) and on variance bounds (Section~\ref{sec:var}). In each case, we indicate where the identity function plays its role and how, by replacing it with $(\log p(x))'$ for some general density $p$, the characterization that seemed tailor-made for the normal distribution can actually be extended to other distributions. Bearing in mind that these extensions are actually based on the location score function, we then explain in Section~\ref{sec:scale} how alternative characterizations can be obtained by rather looking at the scale score function. We conclude the paper with final comments in Section~\ref{sec:final}.

\section{FOUR CHARACTERIZATIONS FROM DIFFERENT RESEARCH DOMAINS}

\subsection{Maximum likelihood characterization}\label{sec:MLE}

We call \emph{location MLE characterization} the characterization of a  probability distribution via the structure of the Maximum Likelihood Estimator (MLE) of the location parameter.   \cite{gauss1809theoria} showed that,  in a location family $p(x -\mu)$ with differentiable density $p$, the MLE for $\mu$ is the sample mean $\bar{x}=\frac{1}{n}\sum_{i=1}^nx_i$ for all samples $(x_1,\ldots,x_n)$ of all sample sizes $n$ if, and only if, $p$ is the normal density. This result has been successively refined in two directions. On the one hand,  several authors have worked towards weakening the regularity assumptions on the class of densities $p$ considered; for instance Gauss requires differentiability  while \cite{teicher1961maximum} only requires continuity. On the other hand, many authors have  lowered the sample size necessary for the characterization to hold, in other words, the ``always''-statement from the Introduction. For instance, Gauss requires that the sample mean be MLE for the location parameter for all sample sizes simultaneously, \cite{teicher1961maximum} only needs that it be MLE for samples of sizes 2 and 3 at the same time, while \cite{azzalini2007gauss} only require that it be so for a single fixed sample size $n \geq 3$. We refer the reader to Section 1.1 of \cite{duerinckx2014maximum} for a complete list of references on the topic.

Location MLE characterizations for other distributions (Laplace, Gumbel, among others) were first proposed on a case-by-case ad hoc basis, before \cite{duerinckx2014maximum} unified all existing results and established the most general MLE characterization results. We shall here not delve into the technical details of \cite{duerinckx2014maximum}, but rather take a look at the proof by \cite{azzalini2007gauss} and explain how the identification of the identity function and a subsequent replacement with $\varphi_p(x)=(\log p(x))'$ directly leads to a generalization of location MLE characterizations.

Letting $g(x-\mu)$ denote a density over $\R$ with location parameter $\mu\in\R$, and assuming differentiability of $g$, the starting point of the proof by \cite{azzalini2007gauss} consists in considering the score equation 
\begin{equation}\label{AGeq}
\sum_{i=1}^n\varphi_g(x_i-\bar{x})=0
\end{equation}
 with $\varphi_g(x)=(\log g(x))'$ for all samples $x_1,\ldots,x_n$ of a fixed sample size $n\geq3$. The particular choices $x_1=x_2=\cdots=x_n=0$, $x_1=a=-x_2,x_3=x_4=\cdots=x_n=0$ and $x_1=a,x_2=b,x_3=-a-b,x_4=\cdots=x_n=0$ with $a,b\in\R$  lead to the functional Cauchy equation $\varphi_g(a+b)=\varphi_g(a)+\varphi_g(b)$ for all $a,b\in\R$. Its unique solution is $\varphi_g(x)= cx$ for some real constant $c$, which is precisely the identity function characterization of the normal distribution that we mentioned in the Introduction, up to the constant $c$.

Now, how can we extend this result to a general density $p$? The key idea lies in the fact that the score equation~\eqref{AGeq} actually is a system of two equations, which in terms of a general density $p$ reads
\begin{equation}\label{oureq}
\sum_{i=1}^n\varphi_g(x_i-\mu)=0\quad\mbox{subject to}\quad \sum_{i=1}^n\varphi_p(x_i-\mu)=0.
\end{equation}
The second equation was somewhat hidden in~\eqref{AGeq} under the form $\sum_{i=1}^n\varphi_g(x_i-\mu)=0\quad\mbox{subject to}\quad \sum_{i=1}^n(x_i-\mu)=0.$ Writing $\alpha_i=\varphi_p(x_i-\mu)$ and assuming $\varphi_p$ to be monotone with image $\R$,   equations~\eqref{oureq} can be rewritten as
\begin{equation}\label{AGnew}
\sum_{i=1}^n\varphi_g\circ\varphi_p^{-1}(\alpha_i)=0\quad\mbox{subject to}\quad \sum_{i=1}^n\alpha_i=0.
\end{equation}
Comparing~\eqref{AGnew} with~\eqref{AGeq}, we notice that both are actually the same set of equations thanks to the monotonicity assumption together with the fact that the $\alpha_i$ span over $\R$. Consequently, we find that $\varphi_g\circ\varphi_p^{-1}(x)=cx$, leading to $g$ equal to $p^c$ with $c$ necessarily positive and hence a location MLE characterization for $p$. From the similarity with the normal proof, it inherits validity for all samples from a fixed sample size $n\geq3$. 

We conclude this section with a few comments. The monotonicity assumption and the fact that $\varphi_p$ maps $\R$ onto all $\R$ are natural extensions from the normal log-density being the identity. Actually, the second requirement may be weakened by only asking that $\varphi_p$ crosses the $x$-axis (otherwise the equation $\sum_{i=1}^n\varphi_p(x_i-\mu)=0$ would have no solution). We refer the reader to \cite{duerinckx2014maximum} for a formal proof of the latter statement. Strict monotonicity and crossing the $x$-axis are actually two requirements that define the class of \emph{strong unimodal} or \emph{log-concave} densities. The location MLE characterization of the normal distribution, initiated by the father of the  Gaussian law himself, thus can more or less readily be extended to this broader family of distributions thanks to the understanding of the role played by the ``normal'' identity function.

\subsection{Singularity of the Fisher information matrix in skew-symmetric distributions}\label{sec:info}

Nice as it is, the symmetric shape of the normal distribution also has its drawbacks, as it does not allow modelling data exhibiting skewness. Consequently, many proposals of transformed normal distributions have been brought forward in the literature. One of the most famous is the \emph{skew-normal} of \cite{azzalini1985class} with density
\begin{equation}\label{sn}
x\mapsto 2 \phi_{\mu,\sigma}(x) \Phi_{0,1}\left(\delta \frac{(x-\mu)}{\sigma}\right),\quad x\in\R,
\end{equation}
where $\Phi_{\mu,\sigma}$ stands for the cumulative distribution function (cdf) associated with the normal density $\phi_{\mu,\sigma}$ and $\delta\in\R$ plays the role of a skewness parameter. At $\delta=0$ we retrieve the normal distribution, and all non-zero values of $\delta$ lead to a skewed distribution. Many further papers have studied various aspects of the skew-normal, and generalizations to other so-called skew-symmetric densities have been proposed both in the univariate and multivariate settings. Scalar skew-symmetric densities are of the form 
$$
x\mapsto \frac{2}{\sigma} p\left(\frac{x-\mu}{\sigma}\right) \Pi\left(\frac{x-\mu}{\sigma},\delta\right),\quad x\in\R,
$$ 
where $p$ is a symmetric density to be skewed and the skewing function $\Pi$ satisfies $\Pi(z,\delta)+\Pi(-z,\delta)=1\,\forall z,\delta\in\R$ and $\Pi(z,0)=\frac{1}{2}\,\forall z\in\R$. The most typical choice of skewing function is $F\left(\delta\frac{(x-\mu)}{\sigma}\right)$ for some symmetric univariate cdf $F$, see~\eqref{sn} where $F=\Phi_{0,1}$. We refer the interested reader to \cite{azzalini2014skew} for a recent overview on skew-symmetric distributions. 

Besides its nice stochastic properties, the skew-normal is also infamous for an inferential peculiarity. In the vicinity of symmetry, that is, when $\delta=0$, the Fisher information matrix associated with the model~\eqref{sn} is singular with rank 2 instead of 3, due to a collinearity between the scores for location and skewness. Straightforward manipulations show that both these scores are proportional to (the identity function) $\frac{x-\mu}{\sigma}$. This singularity prevents for instance the construction of the likelihood ratio test for normality against skew-normality.  \cite{azzalini1985class} proposed a reparameterization that avoids this issue in the scalar case, \cite{arellano2008centred} extended this idea to the multivariate setting, and  \cite{hallin2014skew} have suggested an alternative reparameterization in the scalar case. The unpleasant Fisher information singularity issue, and the ensuing difficulty of building efficient tests for normality, has received a lot of attention in the literature. Mentioned, from the very beginning, by \cite{azzalini1985class}, it is discussed, in the univariate and multivariate skew-normal context, by \cite{azzalini1999statistical}, \cite{pewsey2000problems}, \cite{chiogna2005note}, \cite{ley2010singularity} and \cite{hallin2012skew} among  others. 

A long time  open question in the literature was which other skew-symmetric distributions would suffer from this type of singularity. For instance, \cite{azzalini2003distributions} have proposed the skew-$t$ distribution by using the skewing function $T_{\nu+1}\left(\delta\frac{(x-\mu)}{\sigma}\frac{\nu+1}{\nu+\sigma^{-2}(x-\mu)^2}\right)$ with $T_{\nu+1}$ the cdf of the Student distribution with $\nu+1>1$ degrees of freedom, and noticed ``\textit{It was a pleasant surprise to find that in the present setting the behaviour of the log-likelihood function was to be much more regular, at least for those numerical cases which we have explored}''. The alerted reader will by now have discovered what conveys this seemingly special ``property'' to the skew-normal (and hence the normal): the presence of the identity function inside $\Phi_{0,1}\left(\delta \frac{(x-\mu)}{\sigma}\right)$. Indeed, location and skewness scores at $\delta=0$ in the skew-normal case respectively are given by $\frac{x-\mu}{\sigma^2}$ and $\sqrt{2/\pi}\frac{x-\mu}{\sigma}$. Now, starting from a symmetric density $p$, its location score function is $-\sigma^{-1}\varphi_p\left(\frac{x-\mu}{\sigma}\right)$ and, consequently, a skew-$p$ density will be singular if it is of the form
\begin{equation}\label{skewp}
x\mapsto \frac{2}{\sigma} p\left(\frac{x-\mu}{\sigma}\right) F\left(\delta\varphi_p\left(\frac{x-\mu}{\sigma}\right)\right),\quad x\in\R,
\end{equation}
where the choice of $F$ does not matter. When presented under the form~\eqref{skewp}, the information singularity does not look surprising, and one would expect that, for instance, $\frac{2}{\sigma} q\left(\frac{x-\mu}{\sigma}\right) F\left(\delta\varphi_p\left(\frac{x-\mu}{\sigma}\right)\right)$ for some symmetric density $q$ (not proportional to $p^c$ for some $c>0$ such that $p^c$ is integrable) will not lead to singularity issues. In the skew-normal case, this fact was hidden behind the seemingly innocuous identity function. The latter indeed allows to characterize the skew-normal as the only skew-symmetric distribution suffering from a Fisher information singularity when using the popular skewing function  $F\left(\delta\left(\frac{x-\mu}{\sigma}\right)\right)$  irrespective of the choice of $F$, while such a characterization readily extends to density $p$ if we were to use $F\left(\delta\varphi_p\left(\frac{x-\mu}{\sigma}\right)\right)$. The awareness of the link between Gaussian density and identity function, which we wish to underline in the present paper, allows a direct and complete understanding of the problem.\footnote{For the multivariate case, which is more complex but follows the same reasoning based on identity functions, we refer the interested reader to \cite{hallin2012skew} for a complete solution.}

\subsection{Stein characterizations and Stein's density approach}\label{sec:Stein}

Stein characterizations are an important building block of the famous Stein Method. The  goal of this method, initiated by Charles Stein in 1972 \citep{stein1972bound}, is to provide quantitative
 assessments in distributional {comparison} statements of the form
 $W \approx Z$ where $Z$ follows a known and well-understood
 probability distribution (typically normal as in the Central Limit Theorem) and $W$ is the object
 of interest. In a nutshell, Stein's method consists of  two distinct components, namely
\begin{itemize}
\item[]  \underline{Part~A}: a framework allowing to transform the problem of bounding the 
  error in the approximation of $W$ by $Z$ into a problem of bounding the
  expectation of a certain functional of $W$.  

\medskip

\item[] \underline{Part~B}: a bunch of techniques to bound the
  expectation appearing in Part A; the details of these techniques heavily depend  on the form
  of the functional and on the properties of $W$.
 \end{itemize}

Part B is not of interest for the purpose of this paper, hence we shall not further discuss it and rather refer the interested reader to \cite{ley2017stein} and \cite{ross2011fundamentals}. Part~A directly involves the Stein characterization. For a suitable
operator $\mathcal{A}_Z$  and
for a wide class of functions $\FFF(\AAA_Z)$, the equivalence 
\begin{equation*}
  \label{steingen} W\stackrel{d}{=}Z  \mbox{ if and only if }\E[ \AAA_Z
  f(W)] = 0 \mbox{ for all }f \in \FFF(\AAA_Z),
\end{equation*}
where $\stackrel{d}{=}$ means equality in distribution, represents the Stein characterization of $Z$. The usefulness of such a characterization can readily be seen by noticing that $|  \E[ \AAA_Z
  f(W)]|$ can be used as measure of distance between $Z$ and $W$, and the operator $\AAA_Z  f(\cdot)$ is the abovementioned functional of $W$. \cite{stein1972bound} tackled the normal approximation problem, meaning that $Z$ follows a standard normal distribution, and proposed as operator $\AAA_Z f(x)= f'(x)-x f(x)$. Taking the class $\FFF(\AAA_Z)$ so as to ensure the integrability conditions, one readily sees that
  $$
  \E[ f'(Z)-Z f(Z)]=\int_{-\infty}^{\infty}f'(z)\frac{1}{\sqrt{2\pi}}\exp\left(-\frac{z^2}{2}\right)dz-\int_{-\infty}^{\infty}zf(z)\frac{1}{\sqrt{2\pi}}\exp\left(-\frac{z^2}{2}\right)dz=0
  $$
by integration by parts of the second term. This establishes the sufficiency part of the characterization, and we spare the reader the details of the necessity part. The key behind this integration by parts lies in the fact that $-z \exp\left(-\frac{z^2}{2}\right)$ integrates to $\exp\left(-\frac{z^2}{2}\right)$ which, of course, is due to the fact that $-z$ is the derivative of the log-density of the standard normal. 

The latter observation allows us directly to deduce the form of an operator that should lead to a Stein characterization for a given target density $p$. Replacing $-z$ with $\varphi_p(z)=\frac{p'(z)}{p(z)}$ and letting $Z$ follow the distribution determined by $p$, we readily see that  
  $$
  \E[ f'(Z)+\varphi_p(Z) f(Z)]=\int_{-\infty}^{\infty}f'(z)p(z)dz+\int_{-\infty}^{\infty}\frac{p'(z)}{p(z)}f(z)p(z)dz=0
  $$
by integration by parts (assuming the required minimal integrability conditions). Thus a simple observation and manipulation  leads us to postulate that 
\begin{equation*}
  \label{steingenp} W\stackrel{d}{=}Z\sim p  \mbox{ if and only if }\E[   f'(W)+\varphi_p(W)f(W)] = 0 \mbox{ for all }f \in \FFF(\AAA_Z)
\end{equation*}
is a Stein characterization for $p$, where the class $\FFF(\AAA_Z)$ is determined by the regularity conditions needed for the proof. And indeed: this equivalence happens to be characterizing for any differentiable density $p$ and has come to knowledge in the literature under the name \emph{Stein's density approach} as proposed in~\cite{stein2004use} and further studied in~\cite{ley2013stein}. For the sake of completion, we now provide the formulation from the latter paper which is applicable to a very wide class of distributions.

\begin{theorem}[\cite{ley2013stein}]
Let $Z$ be an absolutely continuous random variable with density $p$ satisfying the following conditions: (i) its support $S_p:=\{z\in\R:p(z)\, \mbox{is positive}\}$ is an interval with  closure $\bar{S}_p=[a,b]$ for some $-\infty\leq a<b\leq\infty$, (ii) it is differentiable at every point in $(a,b)$ with derivative $p'(z)$, and (iii) $\int_{S_p}p(z)dz=1$. Associate to each $p$ the class of functions $\mathcal{F}(p)$ of functions $f:\R\rightarrow\R$ such that the mapping $z\mapsto f(z)$ is differentiable on the interior of $S_p$ and $f(a^+)p(a^+)=f(b^-)p(b^-)=0$. Let $W$ be another absolutely continuous random variable. Then $\E[f'(W)+\varphi_p(W)f(W)]=0$ for all $f\in\mathcal{F}(p)$ if and only if either ${\rm P}(W\in S_p)=0$ or ${\rm P}(W\in S_p)>0$ and ${\rm P}(W\leq z\,|\, W\in S_p)={\rm P}(Z\leq z)$ for all $z\in S_p$.
\end{theorem}

We refer to~\cite{ley2013stein} for the proof. We attract the reader's attention to a slight difference between the statement above and the original statement in~\cite{ley2013stein}: the authors there only require  that the mapping $z\mapsto f(z)p(z)$ be differentiable on the interior of $S_p$, and use as operator $\frac{(f(z)p(z))'}{p(z)}$. Our slightly more stringent condition that $f$ be differentiable on $S_p$ allows rewriting this expression of the operator under the form used above.

\subsection{On a result by \cite{cacoullos1982upper} regarding variance bounds}\label{sec:var}

A famous result of \cite{chernoff1981note} states that, if $X\sim \mathcal{N}(0,1)$, then  the inequality 
$$
{\rm Var}[g(X)]\leq {\rm E}[(g'(X))^2]
$$
holds for an absolutely continuous real-valued function $g$ for which  $g(X)$ has finite variance, and the inequality becomes an equality if and only if $g(x)=ax+b$ for some real constants $a$ and $b$. This type of inequality falls under the umbrella of variance bounds and is useful for solving variations of the classical isoperimetric problem. This result has stimulated the search for general variance bounds, see for instance \cite{cacoullos1982upper}, \cite{klaassen1985inequality}, \cite{afendras2014strengthened}, \cite{ley2016parametric} and \cite{ERS20}. In particular, \cite{cacoullos1982upper} presented the following lemma as basis for upper variance bounds for several densities.

\begin{lem}[\cite{cacoullos1982upper}]\label{lemcac}
Let $X$ be a continuous random variable with density function~$p(x)$. Let $g$ and $g'$ be real-valued functions on $\R$ such that $g$ is an indefinite integral of $g'$, and ${\rm Var}[g(X)]<\infty$. Then
\begin{equation}\label{lemcacexp}
{\rm Var}[g(X)]\leq \int_0^\infty\int_t^\infty x p(x) [g'(t)]^2dx dt-\int_{-\infty}^0\int_{-\infty}^t xp(x)[g'(t)]^2dx dt.
\end{equation}
\end{lem}
This result is a direct consequence of 
\begin{equation}\label{lemcacexp1}
{\rm Var}[g(X)]={\rm Var}\left[\int_0^Xg'(t)dt\right]\leq{\rm E}\left[\left(\int_0^Xg'(t)dt\right)^2\right]\leq{\rm E}\left[X\int_0^X(g'(t))^2dt\right],
\end{equation}
where the last inequality follows from the Cauchy-Schwarz inequality.  Expression~\eqref{lemcacexp} is then readily obtained by writing out explicitly the integrals and switching the integration order. An equality in Lemma~\ref{lemcac} corresponds to ${\rm E}[g(X)]=g(0)$ (first inequality) and $g'(t)\propto 1$ (second, Cauchy-Schwarz, inequality), in other words occurs if and only is $g$ is linear and  ${\rm E}[X]=0$ under $p$. When $p$ is the standard normal density in~\eqref{lemcacexp},  it readily follows that $\int_t^\infty x p(x)dx=p(t)$ and $-\int_{-\infty}^t x p(x)dx=p(t)$ and consequently Lemma~\ref{lemcac} yields that ${\rm Var}[g(X)]\leq ({\rm E}[(g'(X))^2])$ with equality if and only if $g$ is linear. \cite{cacoullos1982upper} also applies this result to the exponential distribution but, in all generality, the lemma is designed for any continuous density $p$. The upper bound he gets for the exponential distribution is however far from optimal (see Section~\ref{varb2}), and the reason lies in the \emph{a priori} hidden fact that Lemma~\ref{lemcac} is actually designed for the normal distribution. This can be recognized by the presence of the identity function $x$ in the right-hand side integrals in~\eqref{lemcacexp} and by the fact that equality holds when $g$ is a linear (hence nearly identity) function. Only when $p$ is the normal density do we have that  $\int_t^\infty x p(x)dx=p(t)$ which is the key element for obtaining a sharp upper bound. Intuitively, this can be seen as the most direct way to get to the density $p$ and hence the expectation in the upper bound, and any superfluous terms would yield worse bounds.  

Again, the alerted reader will now have noted how to improve on Cacoullos' approach for general densities $p$, namely by replacing $x$ with $-\varphi_p(x)$ in~\eqref{lemcacexp}, which however requires a clever prior replacement in~\eqref{lemcacexp1}. This is achieved as follows:
\begin{eqnarray*}
{\rm Var}[g(X)]={\rm Var}\left[\int_0^{-\varphi_p(X)}(g\circ(-\varphi_p)^{-1})'(t))dt\right]&\leq&{\rm E}\left[-\varphi_p(X)\int_0^{-\varphi_p(X)}((g\circ(-\varphi_p)^{-1})'(t))^2dt\right]\\
&=&{\rm E}\left[-\varphi_p(X)\int_0^{-\varphi_p(X)}\frac{(g'((-\varphi_p)^{-1}(t)))^2}{((-\varphi_p)'((-\varphi_p)^{-1}(t)))^2}dt\right]
\end{eqnarray*}
where the monotonicity of $\varphi_p(X)$ is crucial. An equivalent to Lemma~\ref{lemcac} is then readily written down. The key element in obtaining sharp upper bounds for any density $p$ is that, here, $\int_{t}^{\infty} -\varphi_p(x) p(x)dx=p(t)$, which (after some manipulations) leads to the sharp upper variance bound
\begin{eqnarray*}
{\rm Var}[g(X)]\leq {\rm E}\left[\frac{(g'(X))^2}{(-\varphi_p)'(X)}\right]
\end{eqnarray*}
with equality if and only if $g\propto \varphi_p$. This result (inequality plus equality statement) perfectly extends the famous variance bound for the Gaussian distribution and provides insights as to why this approach by Cacoullos actually worked out for the Gaussian, but not for other distributions. The steps  showcased here are a simplified version of the proof provided in~\cite{ley2016parametric}; we refer the reader to that paper for rigorous conditions and for a discussion on the sharpness of such bounds as compared to competitors from the literature.

\section{FURTHER EXTENSIONS BY MEANS OF THE SCALE SCORE FUNCTION}\label{sec:scale}

So far our extensions of characterization theorems of the normal distribution to any other continuous distribution with density $p$ have been based on the score function $\varphi_p(x)=p'(x)/p(x)$ as natural extension of $-x$, the normal score function. It is important to keep in mind that these are  location score functions, as described in the Introduction, and consequently all described characterizations are location-based. For certain distributions this may lead to somewhat artificial results if, for instance, they do not contain a location parameter under their natural form. Think for example of the negative exponential distribution with density $p(x)=\lambda \exp(-\lambda x)$ over $\R^+$, where $\lambda>0$ is a scale parameter. Since it is only defined over the positive real halfline, it does not contain a location parameter by nature. This is also reflected by its constant score function $\varphi_p(x)=-\lambda$ (or $-1$ if we consider the standardized form). Adding a location parameter is of course possible and leads to $p(x)=\lambda \exp(-\lambda (x-\mu))\mathcal{I}(x\geq\mu)$ with $\mu\in\R$ and $\mathcal{I}$ the indicator function. Under this form one speaks of the two-parameter exponential distribution and the location score function should, in principle, involve a weak derivative or derivative in the sense of distributions that also differentiates the indicator function. It is highly questionable if this is the most natural way to treat the exponential distribution and, more generally, distributions over $\R^+$. Instead, it would rather seem appealing to consider the scale score function $\frac{d}{d\lambda}\log\left(\lambda \exp(-\lambda x)\right)=\frac{1}{\lambda}-x$. For a scale family $\sigma p(\sigma x)$, the scale score (after setting $\sigma=1$) is given by $\psi_p(x)= 1+x \varphi_p(x)$ and substituting this quantity for $\varphi_p$ in the examples considered in the previous sections allows obtaining further interesting extensions. In the rest of this section we shall briefly discuss the four previous topics under the light of scale-based characterizations and general parametric parameterizations.

\subsection{MLE characterizations}

Various papers have provided MLE characterizations with respect to the scale parameter. \cite{teicher1961maximum} shows that, in a scale family $\sigma p(\sigma x)$ and under some regularity assumptions, the MLE for the scale parameter $\sigma$ is the sample mean $\bar{x}$  for all samples $x_1,\ldots,x_n$ over $\R^+$ of all sample sizes $n$ if and only if $p$ is the exponential density, while if it corresponds to the square root of the sample arithmetic mean of squares $\left(\frac{1}{n}\sum_{i=1}^nx_i^2\right)^{1/2}$, then this characterizes the normal distribution over $\R$. \cite{marshall1993maximum} extend the characterization of \cite{teicher1961maximum} from the negative exponential to the Gamma distribution. \cite{duerinckx2014maximum} provide a general characterization result for scale families that incorporates all those from the literature. These authors also provide a general MLE characterization for one-parameter group families of the form $H'_\theta(x) p(H_\theta(x))$ where the parameter of interest $\theta$ can take on diverse roles and $H_\theta$ is a differentiable transformation. 

\subsection{Fisher information singularity issue}

The Fisher information singularity within skew-symmetric distributions has only been studied from what we call a location-based view. This is due to the fact that the notorious singularity in the skew-normal case is due to a collinearity between the scores for location and skewness, and consequently  papers such as \cite{hallin2012skew, hallin2014skew} studied the singularity from this viewpoint. We shall therefore consider  here for the first time a skewness-scale induced singularity. It is easy to see that skew-symmetric densities of the form 
\begin{equation}\label{scass}
\frac{2}{\sigma}q\left(\frac{x-\mu}{\sigma}\right) F\left(\delta\psi_p\left(\frac{x-\mu}{\sigma}\right)\right),\quad x\in\R,
\end{equation}
where $p, q$ are symmetric densities and $F$ is  some univariate symmetric cdf, suffer from a singular Fisher information when $\delta=0$ if and only if the scores for scale and skewness are collinear almost everywhere, i.e., $\psi_q(x)=c_1\psi_p(x)+c_2$ a.e. in $x\in\R$ and for some real constants $c_1,c_2$ (since the location score $\varphi_q(x)$ is an odd function, contrary to the even scale score and, here, the also even skewness score at $\delta=0$). The latter equation can be re-expressed under the form
$$
\varphi_q(x)=\frac{(c_1+c_2-1)}{x}+c_1\varphi_p(x)\,{\rm a.e.}.
$$
The solution to this first-order differential equation is $q(x)= dx^{c_1+c_2-1}p^{c_1}(x)$ a.e. and some normalizing constant $d>0$. Hence, for all values $c_1,c_2\in\R$ for which $x^{c_1+c_2-1} p^{c_1}(x)$ is integrable and symmetric, the density $dx^{c_1+c_2-1}p^{c_1}(x)$ leads to a singular Fisher information in the model~\eqref{scass}. The symmetry requirement reduces the possible values of $c_1+c_2$ to odd integers. For the sake of illustration, when $p$ is the normal density, $x^{c_1+c_2-1} \exp\left(-c_1\frac{x^2}{2}\right)$ is integrable for all~$c_1>0$ and $c_2$ such that $c_1+c_2\geq1$ and odd.

\subsection{Stein characterizations}

For exponential approximation problems, \cite{chatterjee2011exponential} use and combine two different Stein characterizations of the negative exponential density $p(x)=\exp(-x)$ over $\R^+$. The first involves the operator $f'(x)-f(x)$ and is applied, in Part B of Stein's Method, for $x\in[0,1]$ while the second considers $xf'(x)-(x-1)f(x)$ and is applied for $x>1$. For the sake of readability we do not specify the regularity assumptions on $f$ and refer to \cite{chatterjee2011exponential} for that purpose. The reader will have noticed that the $-1$ appearing in the first operator corresponds to the location score function $\varphi_p(x)=-1$ while the second operator is based on $\psi_p(x)=1-x$. Thus, without a proper mention in that paper, the authors actually obtained improved upper bounds for exponential approximation by combining location- and scale-based operators. We will not delve here into deeper structural reasons for the $x$ appearing in $xf'(x)$ in the second characterization, and refer the interested reader to \cite{ley2017stein} for more information about setting up useful operators in Stein's Method.

General parametric Stein characterizations, based on parameters of interest of other natures than location and scale, have been studied in \cite{ley2016parametric} and \cite{ley2016general}. The latter paper also develops a link between typical operators from the literature, such as those from \cite{chatterjee2011exponential}, and the operators obtained by adopting the (till then not considered) parametric viewpoint.

\subsection{Variance bounds}\label{varb2}

General parametric variance bounds have been studied in detail in \cite{ley2016parametric}, where the scale case is given particular attention. Since previously no mention on scale-based variance bounds has been made in the literature, this paper compares the resulting bounds to those of \cite{cacoullos1982upper} and \cite{klaassen1985inequality}, noting in particular that the scale-based bounds clearly improve on the Cacoullos bounds and in many situations on the Klaassen bounds. This  underlines the strength of the parametric approach, whose wealth is further undermined in \cite{ley2016parametric} via novel skewness-based variance bounds. It is notable that \cite{cacoullos1982upper} noticed that his variance bounds for the exponential distribution were not very sharp by having recourse to~\eqref{lemcacexp} (in fact, they reached equality if and only if the function $g$ is constant since ${\rm E}[X]$ is not zero in that case) and proposed a way to lower the bound (yielding equality for $g$ linear). However, this was no structural improvement grounded on the nature of the exponential distribution as a scale-based distribution, and hence cannot match the upper bounds from \cite{ley2016parametric}.

\section{FINAL COMMENTS}\label{sec:final}

We hope to have conveyed through the previous examples from very different topics the important message that many characterizations of the normal distribution and, consequently, the seemingly special role of the normal distribution, are (at least to a large degree) to be attributed to the fact that its score function is the identity function which happens to appear in many circumstances. While a general score function of the form $\frac{p'(x)}{p(x)}$ would immediately hint at a special role played by the density $p$, the same does not hold true for $-x$ unless one is aware that $-x=\frac{\phi_{0,1}'(x)}{\phi_{0,1}(x)}$. Keeping this in mind, many results can be better understood and the theory can move forward more quickly.  There exist several further situations where this observation turns out to be useful, for instance in the definition of a generalized Fisher information distance and ensuing information inequalities (see \cite{ley2013stein}) or in the extension to any target $p$   of the normal characterization provided in \cite{nourdin2009density}, see Theorem 2 of~\cite{kusuoka2012stein}.

\bibliographystyle{plain}
\nocite{*}
\bibliography{bibi}

\end{document}